\documentclass[12pt,oneside]{amsart}
\usepackage{amssymb}
\usepackage{amsfonts}
\usepackage{amsthm}
\usepackage{amsthm}
\usepackage{amscd}

\numberwithin{equation}{section}
\overfullrule=0pt \theoremstyle{plain}
\newtheorem{theorem}{Theorem}[section]

\newtheorem{lemma}[theorem]{Lemma}
\theoremstyle{definition}

\begin{document}
\large
%
%

\newcommand{\M}{\mathcal{M}_{g,N+1}^{(1)}}
\newcommand{\Teich}{\mathcal{T}_{g,N+1}^{(1)}}
\newcommand{\T}{\mathrm{T}}
\newcommand{\corr}{\bf}
\newcommand{\vac}{|0\rangle}
\newcommand{\Ga}{\Gamma}
\newcommand{\new}{\bf}
\newcommand{\define}{\def}
\newcommand{\redefine}{\def}
\newcommand{\Cal}[1]{\mathcal{#1}}
\renewcommand{\frak}[1]{\mathfrak{{#1}}}
\newcommand{\refE}[1]{(\ref{E:#1})}
\newcommand{\refS}[1]{Section~\ref{S:#1}}
\newcommand{\refSS}[1]{Section~\ref{SS:#1}}
\newcommand{\refT}[1]{Theorem~\ref{T:#1}}
\newcommand{\refO}[1]{Observation~\ref{O:#1}}
\newcommand{\refP}[1]{Proposition~\ref{P:#1}}
\newcommand{\refD}[1]{Definition~\ref{D:#1}}
\newcommand{\refC}[1]{Corollary~\ref{C:#1}}
\newcommand{\refL}[1]{Lemma~\ref{L:#1}}
\newcommand{\refR}[1]{Remark~\ref{R:#1}}
\newcommand{\R}{\ensuremath{\mathbb{R}}}
\newcommand{\C}{\ensuremath{\mathbb{C}}}
\newcommand{\N}{\ensuremath{\mathbb{N}}}
\newcommand{\Q}{\ensuremath{\mathbb{Q}}}
\renewcommand{\P}{\ensuremath{\mathbb{P}}}
\newcommand{\Z}{\ensuremath{\mathbb{Z}}}
\newcommand{\kv}{{k^{\vee}}}
\renewcommand{\l}{\lambda}
\newcommand{\gb}{\overline{\mathfrak{g}}}
\newcommand{\g}{\mathfrak{g}}
\newcommand{\h}{\mathfrak{h}}
\newcommand{\gh}{\widehat{\mathfrak{g}}}
\newcommand{\hh}{\widehat{\mathfrak{h}}}
\newcommand{\ghN}{\widehat{\mathfrak{g}_{(N)}}}
\newcommand{\gbN}{\overline{\mathfrak{g}_{(N)}}}
\newcommand{\tr}{\mathrm{tr}}
\newcommand{\sln}{\mathfrak{sl}(n)}
\newcommand{\gln}{\mathfrak{gl}(n)}
\newcommand{\son}{\mathfrak{so}(n)}
\newcommand{\spn}{\mathfrak{sp}(2n)}

\newcommand{\gl}{\mathfrak{gl}}
\newcommand{\Hwft}{\mathcal{H}_{F,\tau}}
\newcommand{\Hwftm}{\mathcal{H}_{F,\tau}^{(m)}}

\newcommand{\car}{{\mathfrak{h}}}    
\newcommand{\bor}{{\mathfrak{b}}}    
\newcommand{\nil}{{\mathfrak{n}}}    
\newcommand{\vp}{{\varphi}}
\newcommand{\bh}{\widehat{\mathfrak{b}}}  
\newcommand{\bb}{\overline{\mathfrak{b}}}  
\newcommand{\Vh}{\widehat{\mathcal V}}
\newcommand{\KZ}{Kniz\-hnik-Zamo\-lod\-chi\-kov}
\newcommand{\TUY}{Tsuchia, Ueno  and Yamada}
\newcommand{\KN} {Kri\-che\-ver-Novi\-kov}
\newcommand{\pN}{\ensuremath{(P_1,P_2,\ldots,P_N)}}
\newcommand{\xN}{\ensuremath{(\xi_1,\xi_2,\ldots,\xi_N)}}
\newcommand{\lN}{\ensuremath{(\lambda_1,\lambda_2,\ldots,\lambda_N)}}
\newcommand{\iN}{\ensuremath{1,\ldots, N}}
\newcommand{\iNf}{\ensuremath{1,\ldots, N,\infty}}

\newcommand{\MgN}{\mathcal{M}_{g,N}} 
\newcommand{\MgNeki}{\mathcal{M}_{g,N+1}^{(k,\infty)}} 
\newcommand{\MgNeei}{\mathcal{M}_{g,N+1}^{(1,\infty)}} 
\newcommand{\MgNekp}{\mathcal{M}_{g,N+1}^{(k,p)}} 
\newcommand{\MgNk}{\mathcal{M}_{g,N}^{(k)}} 
\newcommand{\MgNekpp}{\mathcal{M}_{g,N+1}^{(k,p')}} 
\newcommand{\MgNekkpp}{\mathcal{M}_{g,N+1}^{(k',p')}} 
\newcommand{\MgNezp}{\mathcal{M}_{g,N+1}^{(0,p)}} 
\newcommand{\MgNeep}{\mathcal{M}_{g,N+1}^{(1,p)}} 
\newcommand{\MgNeee}{\mathcal{M}_{g,N+1}^{(1,1)}} 
\newcommand{\MgNezz}{\mathcal{M}_{g,N+1}^{(0,0)}} 
\newcommand{\MgNi}{\mathcal{M}_{g,N}^{\infty}} 
\newcommand{\MgNe}{\mathcal{M}_{g,N+1}} 
\newcommand{\MgNep}{\mathcal{M}_{g,N+1}^{(1)}} 
\newcommand{\MgNp}{\mathcal{M}_{g,N}^{(1)}} 
\newcommand{\Mgep}{\mathcal{M}_{g,1}^{(p)}} 
\newcommand{\mpt}{\Sigma,P_1,P_2,\ldots, P_N,\Pif} 
\newcommand{\mpp}{\Sigma,P_1,P_2,\ldots, P_N} 
\newcommand{\At}{\widetilde{\Cal A}}
\newcommand{\Lt}{\widetilde{\Cal L}}
\newcommand{\Lc}{\overline{\Cal L}}
\newcommand{\Yt}{\widetilde{Y}}
\newcommand{\bt}{\tilde b}
\newcommand{\sinf}{{\widehat{\sigma}}_\infty}
\newcommand{\St}{\widetilde{S}}
\newcommand{\uni}{\mathcal{U}}
\newcommand{\can}{\mathcal{K}}
\newcommand{\Wh}{\widehat{W}}
\newcommand{\Wt}{\widetilde{W}}
\newcommand{\MegN}{\mathcal{M}_{g,N+1}^{(1)}} 

\newcommand{\Pif} {P_{\infty}}
\newcommand{\Pinf} {P_{\infty}}
\newcommand{\PN}{\ensuremath{\{P_1,P_2,\ldots,P_N\}}}
\newcommand{\PNi}{\ensuremath{\{P_1,P_2,\ldots,P_N,P_\infty\}}}
\newcommand{\Fln}[1][n]{F_{#1}^\lambda}
\newcommand{\tang}{\mathrm{T}}
\newcommand{\Kl}[1][\lambda]{\can^{#1}}
\newcommand{\A}{\mathcal{A}}
\newcommand{\F}{\mathcal{F}}
\newcommand{\U}{\mathcal{U}}
\newcommand{\V}{\mathcal{V}}
\renewcommand{\O}{\mathcal{O}}
\newcommand{\Ae}{\widehat{\mathcal{A}}}
\newcommand{\Ah}{\widehat{\mathcal{A}}}
\newcommand{\La}{\mathcal{L}}
\newcommand{\Le}{\widehat{\mathcal{L}}}
\newcommand{\Lh}{\widehat{\mathcal{L}}}
\newcommand{\eh}{\widehat{e}}
\newcommand{\Da}{\mathcal{D}}
\newcommand{\kndual}[2]{\langle #1,#2\rangle}
\newcommand{\cins}{\frac 1{2\pi\mathrm{i}}\int_{C_S}}
\newcommand{\cinsl}{\frac 1{24\pi\mathrm{i}}\int_{C_S}}
\newcommand{\cinc}[1]{\frac 1{2\pi\mathrm{i}}\int_{#1}}
\newcommand{\cintl}[1]{\frac 1{24\pi\mathrm{i}}\int_{#1 }}
\newcommand{\w}{\omega}
\newcommand{\ord}{\operatorname{ord}}
\newcommand{\res}{\operatorname{res}}
\newcommand{\nord}[1]{:\mkern-5mu{#1}\mkern-5mu:}
\newcommand{\Fn}[1][\lambda]{\mathcal{F}^{#1}}
\newcommand{\Fl}[1][\lambda]{\mathcal{F}^{#1}}
\renewcommand{\Re}{\mathrm{Re}}

\renewcommand{\k}{{k}}
\newcommand{\ce}{{c}}
\newcommand{\ka}{{k}}
\newcommand{\npos}{{p}}

\define\ldot{\hskip 1pt.\hskip 1pt}
\define\ifft{\qquad\text{if and only if}\qquad}
\define\a{\alpha}
\redefine\d{\delta}
\define\w{\omega}
\define\ep{\epsilon}
\redefine\b{\beta} \redefine\t{\tau} \redefine\i{{\,\mathrm{i}}\,}
\define\ga{\gamma}
\define\cint #1{\frac 1{2\pi\i}\int_{C_{#1}}}
\define\cintta{\frac 1{2\pi\i}\int_{C_{\tau}}}
\define\cintt{\frac 1{2\pi\i}\oint_{C}}
\define\cinttp{\frac 1{2\pi\i}\int_{C_{\tau'}}}
\define\cinto{\frac 1{2\pi\i}\int_{C_{0}}}
\define\cinttt{\frac 1{24\pi\i}\int_C}
\define\cintd{\frac 1{(2\pi \i)^2}\iint\limits_{C_{\tau}\,C_{\tau'}}}
\define\cintdr{\frac 1{(2\pi \i)^3}\int_{C_{\tau}}\int_{C_{\tau'}}
\int_{C_{\tau''}}}
\define\im{\operatorname{Im}}
\define\re{\operatorname{Re}}
\define\res{\operatorname{res}}
\redefine\deg{\operatornamewithlimits{deg}}
\define\ord{\operatorname{ord}}
\define\rank{\operatorname{rank}}
\define\fpz{\frac {d }{dz}}
\define\dzl{\,{dz}^\l}
\define\pfz#1{\frac {d#1}{dz}}

\define\K{\Cal K}
\define\U{\Cal U}
\redefine\O{\Cal O}
\define\He{\text{\rm H}^1}
\redefine\H{{\mathrm{H}}}
\define\Ho{\text{\rm H}^0}
\define\A{\Cal A}
\define\Do{\Cal D^{1}}
\define\Dh{\widehat{\mathcal{D}}^{1}}
\redefine\L{\Cal L} \redefine\D{\Cal D^{1}}
\define\KN {Kri\-che\-ver-Novi\-kov}
\define\Pif {{P_{\infty}}}
\define\Uif {{U_{\infty}}}
\define\Uifs {{U_{\infty}^*}}
\define\KM {Kac-Moody}
\define\Fln{\Cal F^\lambda_n}
\define\gb{\overline{\mathfrak{ g}}}
\define\G{\overline{\mathfrak{ g}}}
\define\Gb{\overline{\mathfrak{ g}}}
\redefine\g{\mathfrak{ g}}
\define\Gh{\widehat{\mathfrak{ g}}}
\define\gh{\widehat{\mathfrak{ g}}}
\define\Ah{\widehat{\Cal A}}
\define\Lh{\widehat{\Cal L}}
\define\Ugh{\Cal U(\Gh)}
\define\Xh{\hat X}
\define\Tld{...}
\define\iN{i=1,\ldots,N}
\define\iNi{i=1,\ldots,N,\infty}
\define\pN{p=1,\ldots,N}
\define\pNi{p=1,\ldots,N,\infty}
\define\de{\delta}

\define\kndual#1#2{\langle #1,#2\rangle}
\define \nord #1{:\mkern-5mu{#1}\mkern-5mu:}
\define \sinf{{\widehat{\sigma}}_\infty}
\define\Wt{\widetilde{W}}
\define\St{\widetilde{S}}
\define\Wn{W^{(1)}}
\define\Wtn{\widetilde{W}^{(1)}}
\define\btn{\tilde b^{(1)}}
\define\bt{\tilde b}
\define\bn{b^{(1)}}
%
\define\eps{\varepsilon}    
\define\doint{({\frac 1{2\pi\i}})^2\oint\limits _{C_0}
       \oint\limits _{C_0}}                            
\define\noint{ {\frac 1{2\pi\i}} \oint}   
\define \fh{{\frak h}}     
\define \fg{{\frak g}}     
\define \GKN{{\Cal G}}   
\define \gaff{{\hat\frak g}}   
\define\V{\Cal V}
\define \ms{{\Cal M}_{g,N}} 
\define \mse{{\Cal M}_{g,N+1}} 
\define \tOmega{\Tilde\Omega}
\define \tw{\Tilde\omega}
\define \hw{\hat\omega}
\define \s{\sigma}
\define \car{{\frak h}}    
\define \bor{{\frak b}}    
\define \nil{{\frak n}}    
\define \vp{{\varphi}}
\define\bh{\widehat{\frak b}}  
\define\bb{\overline{\frak b}}  
\define\Vh{\widehat V}
\define\KZ{Knizhnik-Zamolodchikov}
\define\ai{{\alpha(i)}}
\define\ak{{\alpha(k)}}
\define\aj{{\alpha(j)}}

\newcommand{\wf}[1]{\addtocounter{corr}{1}
            \footnote[\arabic{corr}]{#1}
            \marginpar{\footnotesize{** \arabic{corr}\ **}}
            }

\vspace*{-1cm}
\vspace*{2cm}

\title[ Lax operator algebras]{Lax operator algebras}
\author[I.M. Krichever]{Igor M. Krichever}
\address[Igor M. Krichever]{Landau Institute for Theoretical
Physics, ul.Kosygina 2, Moscow, 119334, Russia,  and Columbia
University, 2960 Broadway, New York, NY 10027-6902, USA}
\email{krichev@math.columbia.edu}
\author[O.K. Sheinman]{Oleg K. Sheinman}
\thanks{Research is supported in part by National Science Foundation
under the grant DMS-04-05519 (I.M.Krichever), by Russian Foundation
for Basic Research under the grant 05-01-00170, and by the program
"Mathematical Methods in Nonlinear Dynamics" of the Russian Academy
of Science (O.K.Sheinman)}
\address[ Oleg K. Sheinman]{Steklov Mathematical Institute, ul. Gubkina, 8,
Moscow, 119991 Russia and Independent University of Moscow, Bolshoi
Vlasievskii per. 11, Moscow, Russia} \email{sheinman@mi.ras.ru}

\begin{abstract}
In this paper we develop a general concept of Lax operators on
algebraic curves introduced in \cite{Kr_Lax}. We observe that the
space of Lax operators is closed with respect to their usual
multiplication as matrix-valued functions. We construct the
orthogonal and symplectic analogs of Lax operators, prove that they
constitute almost graded Lie algebras and construct local central
extensions of those Lie algebras.

\end{abstract}

\maketitle

\section{Introduction}\label{S:Intr}

A general concept of Lax operators on algebraic curves was
proposed by one of the authors in \cite{Kr_Lax}, where the theory
of conventional Lax and zero curvature representations with the
{\it rational spectral parameter} was generalized on the case of
algebraic curves $\Gamma$ of arbitrary genus $g$. The linear space
of such operators associated with an effective divisor $D=\sum_k
n_kP_k, \ P_k\in \Gamma$ was defined as the space of meromorphic
$(n\times n)$ matrix-valued functions on $\Gamma$ having poles of
multiplicity at most $n_k$ at the points $P_k$, and at most simple
poles at $ng$ points $\gamma_s$ more. The coefficients of the
Laurent expansion of those matrix-valued functions in the
neighborhood of a point $\gamma_s$ had to obey certain linear
constrains parametrized by the point $\alpha_s$ of a certain
projective space (see relations \refE{i0}-\refE{ii} below).


According to \cite{Tyur}, the generic sets $(\ga_s,\a_s)$
parameterize stable rank $n$ degree $ng$ framed holomorphic vector
bundles $B$ on $\Gamma$. In \cite{Kr_Lax} it was noticed that the
requirements on the form of Lax operators at the points $\ga_s$
mean that these operators can be seen as meromorphic sections of
the bundle $End(B)$ with $D$ as the divisor of poles. As an easy
consequence of this remark, we have that Lax operators having
poles of arbitrary orders at the points $P_k$ form an algebra with
respect to the usual pointwise multiplication.

In the simplest case of two marked points, $D=P_++P_-$, this
enables us to define in the algebra of the corresponding operators
an almost graded structure generalizing the graded structure of
the classic affine algebra $\widehat{\gl(n)}$. Recall that a Lie
algebra $\V$ is called  {\it almost graded} if $\V=\oplus\V_i$
where $\dim\,\V_i<\infty$ and
$[\V_i,\V_j]\subseteq\oplus_{k=i+j-k_0}^{k=i+j+k_1}\V_k$ where
$k_0$ and $k_1$ do not depend on $i$, $j$.

A general notion of almost graded algebras and modules over them
was introduced in \cite{KNFa}-\cite{KNFc} where the
generalizations of Heisenberg  and Virasoro algebras were
introduced. In the series of works, see \cite{ShN65} for the
review, the almost graded analogues of classic affine Lie
algebras, called Krichever-Novikov (KN) current algebras, were
investigated. It is natural to consider the algebra of Lax
operators having poles at two points as a generalization of the KN
$\gl(n)$-algebra.

A central extension of $\V$ is called {\it local} if it is an
almost graded Lie algebra itself. Local central extensions are
given by local 2-cocycles. A 2-cocycle $\ga$ is called local if a
$K\in\Z$ exists such that $\ga(\V_i,\V_j)=0$ for $|i+j|>K$. The
notion of local cocycle is introduced in \cite{KNFa}. In the same
place its uniqueness was conjectured and the proof for the
Virasoro-type algebras was outlined; a complete proof is given in
\cite{SchlCocycl,SchlMMJ}. The locality condition is important in
consideration of analogues of highest-weight representations.

To construct the orthogonal and symplectic analogues of Lax
operators is the main goal of the present paper. In these cases,
the Lax operators do not form any associative algebra, they only
form a Lie algebra. For all classic Lie algebras $\g$ the
corresponding Lax operator algebras can be considered as a
"twisted" version of the Krichever-Novikov current algebras and
loop algebras.

In \refS{curr} we prove the multiplicative properties of
$\gl(n)$-valued Lax operators, introduce $\g$-valued Lax operators
for $\g=\frak{so}(n)$ and $\g=\frak{sp}(2n)$ and prove that they
are closed with respect to the pointwise bracket.

In \refS{almgrad}, we establish the almost graded structure on Lax
operator algebras, and show that $\dim\,\V_i=\dim\,\g$, as for KN
algebras.

In \refS{ext}, on every type of a Lax operator algebra, we define
a 2-cocycle and prove its locality.

Authors are thankful to M.Schlichenmaier for his fruitful
criticism.

\section{Lax operators and their Lie bracket}\label{S:curr}
\subsection{Lax operator algebras for $\frak{gl}(n)$ and $\sln$}\label{SS:gln}
Following \cite{Kr_Lax}, call by a {\it Lax operator} with  the
Tyurin parameters $\{\a_s,\ga_s|s=1,\ldots,gr\}$ a $\gln$-valued
function $L$ on $\Gamma$, holomorphic outside $P_\pm$ and
$\{\ga_s|s=1,\ldots,gr\}$, having at the latters at most simple
poles:
\begin{equation}\label{E:i0}
 L=\frac{L_{s,-1}}{z-z_s}+L_{s0}+O(z-z_s),\quad z_s=z(\ga_s),
\end{equation}
and such that
\newline
(i) $L_{s,-1}=\a_s\b_s^t$ where $\b_s^t$ is a row, and
\begin{equation}\label{E:i}
 \quad \tr\,L_{s,-1}=\b_s^t\a_s=0,
\end{equation}
where $\a_s\in\C^n$ is fixed, $\b_s\in\C^n$ is arbitrary, and the
upper $t$ denotes a transposed matrix. In particular, $L_{s,-1}$
has rank $1$;
\newline
(ii) $\a_s$ is an eigenvector of the matrix $L_{s0}$:
\begin{equation}\label{E:ii}
  L_{s0}\a_s=\k_s\a_s.
\end{equation}
\begin{lemma}\label{L:Higgs2}
Let $L'$  and $L''$  satisfy the conditions \refE{i0}-\refE{ii}.
Then $L=L'L''$ also satisfies them.
\end{lemma}
\begin{proof} From \refE{i0} we have
\begin{multline}
     L=\frac{L_{s,-1}'L_{s,-1}''}{(z-z_s)^2}+
     \frac{L_{s,-1}'L_{s0}''+L_{s0}'L_{s,-1}''}{(z-z_s)}+\\
 +L_{s,-1}'L_{s1}''+L_{s0}'L_{s0}''+L_{s1}'L_{s,-1}''+O(1).
\end{multline}

From \refE{i} for $L'$  it follows that the first term vanishes:
\[L_{s,-1}'L_{s,-1}''=\a_s({\b_s'}^t\a_s){\b_s''}^t=0. \]
For the second one we have
$L_{s,-1}=L_{s,-1}'L_{s0}''+L_{s0}'L_{s,-1}''=\a_s({\b_s'}^tL_{s0}'')+(L_{s0}'\a_s){\b_s''}^t$.
By  \refE{i} for $L'$  it follows $L_{s0}'\a_s=\k_s'\a_s$, hence
$L_{s,-1}=\a_s\b_s^t$  where
$\b_s^t={\b_s'}^tL_{s0}''+\k_s'{\b_s''}^t$. Further on,
$\tr\,L_{s,-1}=({\b_s'}^tL_{s0}''+\k_s'{\b_s''}^t)\a_s=\k_s''{\b_s'}^t\a_s+\k_s'{\b_s''}^t\a_s=0$.

Consider the expression $L_{s,0}\a_s$, where
$L_{s,0}=L_{s,-1}'L_{s1}''+L_{s0}'L_{s0}''+L_{s1}'L_{s,-1}''$.
From the definition of the Lax operators it follows that
$L_{s,-1}''\a_s=0$ and $L_{s0}'L_{s0}''\a_s=\k_s'\k_s''\a_s$. We
also have $L_{s,-1}'L_{s1}''\a_s=\a_s({\b_s'}^tL_{s1}''\a_s)$.
Hence $\a_s$  is an eigenvector of the matrix $L_{s,0}$ with the
eigenvalue $\k_s={\b_s'}^tL_{s1}''\a_s+\k_s'\k_s''$.
\end{proof}

Since the conditions  \refE{i0},\refE{i} are linear, the Lax
operators constitute an associative algebra, hence the
corresponding Lie algebra. The latter is called {\it Lax operator
algebra}.

If the function $L$, besides \refE{i0}---\refE{ii}, satisfies the
condition $\tr\, L=0$, it is called {\it an $\sln$-valued Lax
operator}. Such Lax operators form a Lie algebra.


\subsection{Lax operator algebras for $\frak{so}(n)$}\label{SS:so_n}
For the elements of this Lie algebra we have $X^t=-X$. We
introduce the matrix-function $L$ taking values in $\son$ by the
same expansion as in the \refSS{gln} but change the requirement
(i) of that section for the reason that there is no rank one
skew-symmetric matrix, and change accordingly the requirement
(ii). We omit the index $s$ for brevity and write down the
expansion \refE{i0} in the form
\begin{equation}\label{E:i0'}
 L=\frac{L_{-1}}{z}+L_{0}+O(z)
\end{equation}
where $L_0,L_1,\ldots$ are skew-symmetric. Instead the condition
(i) of \refSS{gln} we require
\begin{equation}\label{E:i'}
 L_{-1}=\a\b^t -\b\a^t
\end{equation}
where $\a\in\C^n$ is fixed, $\b\in\C^n$ is arbitrary and
\begin{equation}\label{E:j'}
 \a^t\a=\b^t\a(=\a^t\b)=0.
\end{equation}
Likewise \refE{ii} we require
\begin{equation}\label{E:k'}
L_0\a=\k\a
\end{equation}
for some complex number $\k$.

We will prove now the closeness of the space of Lax operators with
respect to the Lie bracket in case of $\g=\frak{so}(n)$ in
question. We stress that there is no structure of associative
algebra in this case.
\begin{lemma}\label{L:mult_so2n}
  The properties \refE{i0'}---\refE{k'}
  are invariant with respect to the Lie bracket.
\end{lemma}
\begin{proof}
 1). First, prove the absence of the term with $z^{-2}$. The corresponding
 coefficient is equal to
\begin{align*}[L_{-1},L_{-1}']&=[\a\b^t-\b\a^t,\a\b'^t-\b'\a^t] \\
&=(\b^t\a)(\a\b'^t-\b'\a^t) -
 (\a^t\a)(\b'\b^t-\b\b'^t) - (\a^t\b')(\a\b^t-\b\a^t).
\end{align*}
It vanishes due to \refE{j'} (applied to both $\b$ and $\b'$). We
stress that for the product $L_{-1}L_{-1}'$ the term with $z^{-2}$
does not vanish.

2). Now calculate the term with $z^{-1}$ in the product $LL'$. The
coefficient is equal to
\begin{align*}
 L_{-1}L_{0}'+L_{0}L_{-1}' & = \a(\b^tL_0')-\b(\a^tL_0')+(L_0\a)\b'^t-
                                            (L_0\b')\a^t\\
     &= \a(\b^tL_0')-\b(-\k'\a^t)+\k\a\b'^t-
                                            (L_0\b')\a^t\\
     &=\a(\b^tL_0'+\k\b'^t)-(L_0\b'-\k'\b)\a^t
\end{align*}
(here we used the relation \refE{k'}). We see that it is not of
the required form \refE{i'}. Now consider the corresponding term
of the commutator of $L$ and $L'$. Then we obtain
\begin{align*} [L,L']_0
 &= \a(\b^tL_0'-\b'^tL_0+\k\b'^t-\k'\b^t)-(L_0\b'-L_0'\b-\k'\b+\k\b')\a^t\\
   &=\a\b''^t-\b''\a^t
\end{align*}
where $\b''^t=\b^tL_0'-\b'^tL_0+\k\b'^t-\k'\b^t$.

It is easy to check that $\b''$ satisfies \refE{j'}.

3) Check the requirement \refE{k'} on the eigenvalue of the degree
$0$ matrix coefficient $(LL')_0=L_{-1}L_1'+L_0L_0'+L_1L_{-1}'$.

By \refE{j'} $L_{-1}'\a=0$, hence for the third term we have
$L_1L_{-1}'\a=0$.

For the first term we have
\[  L_{-1}L'_1\a=(\a\b^t-\b\a^t)L_1'\a=\a(\b^tL_1'\a) -
    \b(\a^tL_1'\a).
\]
The last term on the right hand side of this equation vanishes due
to skew-symmetry of $L'_1$.

Thus  \begin{equation}\label{E:eigen_so} (LL')_0\a=\k''\a\ \
\text{where}\ \  \k''=\b^tL_1'\a+\k\k' .
\end{equation}
\end{proof}

\subsection{Lax operator algebras for $\frak{sp}(2n)$}\label{SS:sp2n_new}
For the elements of the symplectic algebra we have $X^t=-\s
X\s^{-1}$ where $\s$ is a nondegenerate skew-symmetric matrix.

Take the expansion for $L$ in the form
\begin{equation}\label{E:i0''}
 L=\frac{L_{-2}}{z^2}+\frac{L_{-1}}{z}+L_{0}+L_1z+L_2z^2+O(z^3)
\end{equation}
(again, we omit the index $s$ for brevity) where
$L_{-2},L_{-1},L_0,L_1,\ldots$ are symplectic matrices, and
\begin{equation}\label{E:i''}
 L_{-2}=\nu\a\a^t\s,\ \ L_{-1}=(\a\b^t +\b\a^t)\s\quad (\nu\in\C,\ \b\in\C^{2n}).
\end{equation}
Likewise \refE{j'} we require
\begin{equation}\label{E:j''}
 \b^t\s\a=0.
\end{equation}
Notice that $\a^t\s\a=0$ due to the skew-symmetry of the matrix
$\s$.

Further on, we require
\begin{equation}\label{E:k''}
L_0\a=\k\a
\end{equation}
for some complex number $\k$.

We impose a new relation
\begin{equation}\label{E:l''}
  \a^t\s L_1\a=0.
\end{equation}

Prove now that the space of Lax operators is closed with respect
to the Lie bracket in the case $\g=\frak{sp}(2n)$ in question. We
stress that also in this case there is no structure of associative
algebra.
\begin{lemma}\label{L:mult_so2n}
The properties \refE{i0''}---\refE{l''} are invariant with respect
to the Lie bracket.
\end{lemma}
\begin{proof}
Let $L''=[L,L']$.

1) The absence of the terms of orders $-4$ and $-3$ in $z$ in
$L''$ follows from the relations $\b^t\s\a=0$, $\a^t\s\a=0$
\refE{j''} only.

2) For the order $-2$ term we have:
\[
   L''_{-2}=(\nu\,'\k-\nu\k'+\b^t\s\b'){\a\a^t\s},
\]
hence it has the form required by the relation \refE{i''} (here
and below $\nu'$, $\b'$, $L_i'$ have the same meaning for $L'$ as
$\nu$, $\b$, $L_i$ for $L$).

3) The form of the order $-1$ term: a straightforward calculation
using \refE{i0''}, \refE{i''} gives
\begin{align*}
 L''_{-1}=& \a(\nu\cdot\a^t\s L_1'-\nu'\cdot\a^t\s L_1 +\b^t\s L'_0
                -\b'^t\s L_0+\k\b'^t\s-\k'\b^t\s)\\
         &+ (-\nu L_1'\a+\nu' L_1\a- L'_0\b
                +L_0\b'+\k\b'-\k'\b)\a^t\s .
\end{align*}
Denote the second bracket by $\b''$. Then the relations $L_1^t=-\s
L_1\s^{-1}$, $L_0^t=-\s L_0\s^{-1}$ (obeyed by symplectic
matrices) imply that the first bracket is equal to $\b''^t\s$.
Hence\[
  L''_{-1}=(\a\b''^t+\b''\a^t)\s ,
\]
where
\[
   \b''=-\nu L_1'\a+\nu' L_1\a- L'_0\b+L_0\b'+\k\b'-\k'\b .
\]

Prove that $\b''^t\s\a=0$. Making use of the above expression for
$\b''^t\s$ we find
\[
  \b''^t\s \a = \nu\cdot\a^t\s L_1'\a-\nu'\cdot\a^t\s L_1\a +
      \b^t\s L'_0\a-\b'^t\s L_0\a+\k\b'^t\s\a-\k'\b^t\s\a.
\]
The first two terms of the just obtained expression vanish by the
relation \refE{l''} applied to $L$ and $L'$. To the second pair of
terms, we apply the relations $L_0\a=\k\a$ and $L_0'\a=\k'\a$,
after that all remaining terms vanish by the relations \refE{j''}.

4). Prove the relation \refE{k''} on the eigenvalues of the degree
zero term. By definition,
\[
 (LL')_0 = \nu\a\a^t\s L_2'+(\a\b^t+\b\a^t)\s L_1'+L_0L_0'+
 L_1(\a\b'^t+\b'\a^t)\s+ \nu\,'L_2\a\a^t\s.
\]
The last pair of terms obviously vanishes having been multiplied
by $\a$ from the right. We obtain
\[
 (LL')_0\a = \nu\a\a^t\s L_2'\a+\a\b^t\s L_1'\a+\b\a^t\s L_1'\a+\k\k'\a.
\]
The third summand is equal to $0$ by \refE{l''}. Thus $\a$ is the
eigenvector for a zero degree term even in the product $LL'$:
\[
  (LL')_0\a = \a(\nu\cdot\a^t\s L_2'\a+\b^t\s L_1'\a+\k\k'),
\]
and in the commutator we obtain
\begin{equation}\label{E:kappa}
   L''_0\a= \a(\nu\cdot\a^t\s L_2'\a-\nu'\cdot\a^t\s L_2\a +
   \b^t\s L_1'\a-\b'^t\s L_1\a).
\end{equation}

5). Conservation of the relation $\a^t\s L_1\a=0$ in the product
and commutator.

For the product, by definition
\[
  (LL')_0=L_{-2}L_3'+L_{-1}L_2'+L_0L_1'+L_1L_0'+L_2L_{-1}'+L_3L_{-2}'.
\]
Replacing here $L_{-2}$, $L_{-2}'$, $L_{-1}$, $L_{-1}'$ by the
known expressions obtain
\begin{align*}
 \a^t\s (LL')_0\a &= \nu(\a^t\s\a)\a^t\s L_3'\a+
         \left((\a^t\s\a)\b^t+(\a^t\s\b)\a^t\right)\s L_2'\a \\
        & -\k(\a^t\s L_1'\a)+ {\tilde\k}(\a^t\s L_1\a) \\
   & +\a^t\s L_2\left(\a(\b'^t\s\a)+\b'(\a^t\s\a)\right)+
   \nu'\a^t\s L_3\a(\a^t\s\a).
\end{align*}
By the relations \refE{i''}-\refE{l''} this expression vanishes.
\end{proof}


\section{Almost graded structure}\label{S:almgrad}

In this section we consider the following cases: $\g=\sln$,
$\g=\son$, $\g=\spn$, and, to conclude with, $\g=\gln$ which is
somewhat special case with respect to the almost gradeness. A
general definition of the almost graded structure is given in the
introduction. Denote by $\gb$ the Lax operator algebra
corresponding to $\g$.

For all above listed algebras $\g$, except for $\gln$, and
$m\in\Z$ let
\[ \g_m=\{L\in\gb\ |(L)+D\ge 0\},
\]
where $(L)$ is a divisor of the $\g$-valued function $L$, for
$\g=\sln$, $\g=\son$
\[
  D=-mP_+ + (m+g)P_- + \sum_{s=1}^{ng}\ga_s ,
\]
and for $\g=\spn$
\[
  D=-mP_+ + (m+g)P_- + 2\sum_{s=1}^{ng}\ga_s .
\]

We call $\g_m$ {\it the (homogeneous) subspace of degree $m$} of
the Lie algebra $\gb$.

\begin{theorem}\label{T:almgrad} For $\g=\sln ,
\ \mathfrak{so}(n), \ \mathfrak{sp}(2n)$
\begin{align*}
  & 1^\circ . \quad \dim\,\g_m=\dim\,\g.
        \phantom{aaaaaaaaaaaaaaaaaaaaaaaaaaaaaaaaaaa}\\
  & 2^\circ . \quad \gb=\bigoplus\limits_{m=-\infty}^\infty \g_m.\\
  & 3^\circ . \quad [\g_k,\g_l]\subseteq\bigoplus\limits_{m=k+l}^{k+l+
g}\g_m.
\end{align*}
\end{theorem}
\begin{proof}
First prove $1^\circ$. For an arbitrary Lie algebra $\g$, the
Riemann-Roch dimension of all $\g$-valued functions $L$ satisfying
the relation $(L)+D\ge 0$ is equal to $(\dim\g)(ng+1)$. For
$\g=\spn$ it is equal to $(\dim\g)(2ng+1)$. We will prove that for
$\g=\sln$, $\g=\mathfrak{so}(n)$ and an arbitrary $m\in\Z$, there
are exactly $\dim\g$ relations at every pole $\ga_s$, and for
$\g=\mathfrak{sp}(2n)$ the number of those relations is $2\dim\g$.
This will mean that $\dim\g_m=\dim\g$.

Consider first the case $\g=\sln$. The elements of the subspace
$\g_m$ satisfy certain conditions of three kinds coming from
\refE{i}, \refE{ii}: on the residues, eigenvalues and traces of
the matrix $L\in\g_m$, as follows:
\newline\label{conds}
1) at every weak singularity, $L_{-1}=\a\b^t$ which would give
$\dim\g$ relations (since $L_{-1}\in\g$) if the right hand side
was fixed. But the latter depends on the free $n$-dimensional
vector $\b$. Hence we have $\dim\g-n$ conditions at every one of
the $ng$ simple poles $\ga_s$;
\newline
2) at every weak singularity $L_0\a=\k\a$ which gives $n$
conditions. Taking account of one free parameter $\k$ we obtain
$n-1$ condition at every $\ga_s$;
\newline
3) we also have $\tr\, L=0$, i.e. one more relation at every weak
singularity.

Thus we have $(\dim\g-n)+(n-1)+1=\dim\g$ relations at every
$\ga_s$, as required.

For $\g=\mathfrak{so}(n)$ we follow the same line of arguing.
Again, the relation $L_{-1}=\a\b^t-\b\a^t$ \refE{i'} gives
$\dim\g-n$ equations, the relation $L_0\a=\k\a$ \refE{k'} gives
$n-1$ equation, and $\b^t\a=0$ \refE{j'} gives one more equation.
Summing up, we obtain $\dim\g$ equations at every point $\ga_s$.

In the case $\g=\spn$ at every point $\ga_s$ we have the following
conditions:
\begin{align*}
  L_{-2}& =\nu\a\a^t \ \text{:}\quad\quad\quad\quad\ \ \dim\g-1\quad \text{conditions}\
  (\text{one free parameter}\  \nu)\\
  L_{-1}& = (\a\b^t+\b\a^t)\s \ \text{:}\quad \dim\g-2n\ \text{conditions}\
  (2n\ \text{free parameters}\  \b)\\
  L_0\a & =\k\a \ \text{:} \quad\quad\quad\quad\quad\quad
  2n-1\quad\quad \text{conditions}\
  (\text{one free parameter}\  \k)\\
  \b^t\s & \a=0,\ \a^t\s L_1\a=0 \ \text{:}\quad 2\ \text{conditions},
\end{align*}
i.e. $2\dim\g$ conditions at every one of the $ng$ points, as
required.

For $m>0$ and $m<-g$ the subspaces $\g_m$ are linear independent for
the obvious reason: the orders at $P_{\pm}$ of the elements of those subspaces
are different for different $m$.

For $-g\leq m\leq 0$ the linear independence of $\g_m$ follows from
the fact that there are no {\it everywhere holomorphic} Lax
operators. We would like to emphasize that the last argument applies
only to the case of simple Lie algebras. That explains why the case
of the reductive Lie algebra $\gl(n)$ requires some modification
(see below).

The statement $3^\circ$ of the theorem follows from the
consideration of orders at the points $P_\pm$.
\end{proof}

The \refT{almgrad} defines an {\it almost graded structure} on
$\gb$.

\vskip5pt

Consider now the case $\g=\gln$. In this case $\gb$ contains the
subspace of functions taking values in the $1$-dimensional space of
scalar matrices. Let $L$ be such function. By \refE{i} we obtain
$\tr\, L_{-1}=0$. Since $L_{-1}$ is scalar, we obtain $L_{-1}=0$.
Hence, $L$ is holomorphic except at $P_\pm$. Let $\A$ denote the
algebra of meromorphic functions on $\Gamma$ holomorphic except at
$P_\pm$. Then $L\in\A\cdot id$ where $id$ is the unit matrix.
Therefore,
\begin{equation}\label{E:dirs}
{\overline\gln}={\overline\sln}\oplus\A\cdot id .
\end{equation}
In \cite{KNFa}, there was introduced a certain base $\{A_m\}$
(later called the {\it Krichever-Novikov base}) in the space of
such functions. Denote $\A_m=\C A_m$ and set $\gln_m=\sln_m\oplus
(\A_m\cdot id)$. For $m>0$ and $m<-g$ this definition is
equivalent to the above given definition of $\g_m$ (with
$\g=\gln$).

With $\g_m=\gln_m$, the \refT{almgrad} remains to be true for
$\g=\gln$, as it follows from \refE{dirs}. The relation $3^\circ$
of the theorem holds with the different upper limit of summation
determined by the algebra $\A$ (see \cite{KNFa}).

\section{Central extensions of Lax operator algebras}\label{S:ext}

\subsection{Central extensions for the Lax operator algebras over $\mathfrak{gl}(n)$}\label{SS:ext_gln}

The conventional expression for the 2-cocycle giving a central
extension for Krichever-Novikov current algebras (in particular,
loop algebras) is \linebreak $\tr\,\res_{P_+} L\,dL'$. For those
algebras the cocycle is local. A cocycle $\chi$ is called {\it
local} \cite{KNFa} if the constants $\mu'$, $\mu''$ exist such that
$\chi(\g_m,\g_{m'})=0$ unless $\mu'\le m+m'\le\mu''$ where $\g_m$,
$\g_{m'}$ are homogeneous subspaces introduced in the previous
section. In case of Lax operator algebras the above cocycle is by no
means local. In this section we improve it in order to obtain a
local cocycle.

For every operator $L$ the eigenvalue $k$ of its zero degree
component $L_0$ (see \refE{ii}) can be considered as a linear
functional of $L$. We denote this functional by $\k(L)$.
\begin{lemma}\label{L:eig} At every weak singularity, the one-form
$\tr\, L\,dL'$ has at most simple pole and
\begin{equation}\label{E:eig}\res\,\tr\, L\,dL'=\k([L,L']).
\end{equation}
\end{lemma}
\begin{proof} Calculate both parts of the relation explicitly.

1) Using
\[ dL'=-\frac{\a\b'^t}{z^2}+L'_1+\ldots
\]
and \refE{i0} we obtain
\[  L\,dL'=-\frac{\a\b^t\a\b'^t}{z^3}- \frac{L_0\a\b'^t}{z^2} - \frac{L_1\a\b'^t-\a\b^tL_1'}{z} + \ldots\ .
\]
The first term vanishes since $\b^t\a=0$. The second term vanishes
under the symbol of trace, since $L_0\a=\k\a$, and
$\tr\,\a\b'^t=\b'^t\a=0$. The third term gives us the required
residue. We have
\[ \res\,\tr\, L\,dL'= \tr\,(\a\b^tL_1'-L_1\a\b'^t) = \b^tL_1'\a-\b'^tL_1\a.
\]

2) Now calculate the right hand side of the relation in question.
Denote by $[L,L']_0$ the zero degree coefficient of the expansion
\refE{i0} for the commutator $[L,L']$. We have
\begin{align*} [L,L']_0=\ \ &\a\b^tL_1'+L_0L_0'+L_1\a\b'^t \\
             -&\a\b'^tL_1-L_0'L_0-L_1'\a\b^t .
\end{align*}
Multiply both parts of this relation by $\a$ from the right hand
side. Then the third terms in every row vanish since they contain
the factors $\b'^t\a$, $\b^t\a$ which are zeroes by \refE{i}. The
second terms annihilate since $L_0L_0'\a=\k'\k\a$ and
$L_0'L_0\a=\k\k'\a$. Hence
\[ [L,L']_0\a=\a\b^tL_1'\a-\a\b'^tL_1\a = \a(\b^tL_1'\a-\b'^tL_1\a) .
\]
The expression in brackets is a $1\times 1$-matrix, that is a
complex number. Certainly, this complex number is the eigenvalue
of $[L,L']_0$ on the vector $\a$, i.e. $\k([L,L'])$. Its value
exactly meets the above resulting expression for $\tr\,\res\,
L\,dL'$, which completes the proof.
\end{proof}
We want to eliminate the singularities of the 1-form $\tr\, LdL'$
at $\ga_s$ by subtracting another expression for the eigenvalue on
the right hand side of \refE{eig}. Remarkably, we are able to give
this new expression in terms of connections in the holomorphic
vector bundles on $\Gamma$, given in explicit form in
\cite{Kr_Iso}. In what follows, $\L$ denotes the 1-form of such
connection.
\begin{lemma}\label{L:iso}  Let $\mathcal L$ be a matrix-valued 1-form such that locally, at a weak singularity,
\[
 \L=\L_{-1}\frac{dz}{z}+
    \L_0 dz    +\ldots\ ,
\]
 and $\mathcal L$ satisfies the similar requirements as $L$
 (see \refE{i0})-\refE{ii}), with the only one modification:
 ${\tilde\b}^t\a=1$ where $\L_{-1}=\a{\tilde\b}^t$. Then the 1-form
 $\tr\, L\L$ has at most simple pole at $z=0$ and
\[ \res\,\tr\, L\L =\k(L).
\]
\end{lemma}
\begin{proof}
We have
\[ L=\frac{\a\b^t}{z}+L_0+\ldots\ ,\ \
 \L=\left(\frac{\a{\tilde\b}^t}{z}+\L_0+\ldots\right)dz ,
\]
hence
\[ L\L=\left(\frac{\a\b^t\a{\tilde\b}^t}{z^2}+ \frac{\a\b^t\L_0+
       L_0\a{\tilde\b}^t}{z} +\ldots\right)dz .
\]
As above, the first term vanishes. For the second term we have
\[ \res\,\tr\, L\L = \tr\,(\a\b^t\L_0+ L_0\a{\tilde\b}^t) = \b^t\L_0\a + {\tilde\b}^tL_0\a.
\]
Further on, $\L_0\a={\tilde\k}\a$, $L_0\a=\k\a$, hence
$\b^t\L_0\a={\tilde\k}\b^t\a=0$ and
${\tilde\b}^tL_0\a=\k{\tilde\b}^t\a=\k$. This completes the proof.
\end{proof}
\begin{theorem}\label{T:loc}
For every $\L$ satisfying the conditions of the \refL{iso} the
1-form $\tr\, (L\,dL'-[L,L']\L)$  is regular except at the points
$P_\pm$, and the expression
\[  \ga(L,L')= \res_{P_+}\tr\, (L\,dL'-[L,L']\L)
\]
gives a local cocycle on the Lax operator algebra.
\end{theorem}
\begin{proof}
In course of proofs of the \refL{eig} and the \refL{iso}, we have
seen that the 1-forms $\tr\, L\,dL'$ and $\tr\, [L,L']\L$ have
simple poles at every point $\ga_s$, and their residues are equal to
the same quantity $\k_s([L,L'])$. Hence, their difference is regular
at every $\ga_s$.

Assume, at the point $P_+$ we have
\begin{equation}\label{E:expan}
  L(z)=\sum_{i=m}^{\infty} a_iz^i,\
  L'(z)=\sum_{j=m'}^{\infty} b_jz^j,\
  \L(z)=\sum_{k=m_+}^{\infty} c_kz^kdz.
\end{equation}
Then
\[
  L(z)\, dL'(z)=\sum_{p=m+m'}^\infty \left(\sum_{i+j=p}ja_ib_j\right)z^{p-1}dz,
\]
and
\[
  [L(z),L'(z)]\L=\sum_{p=m+m'+m_+}^\infty \left(\sum_{i+j+k=p}[a_i,b_j]c_k\right)z^pdz.
\]
In order one of these 1-forms had a nontrivial residue at the
point $P_+$ it is necessary that either $m+m'\le 0$ or
$m+m'+m_+\le -1$, in other words
\[  m+m'\le \max\,\{0, -1-m_+\}.
\]

If $L$ and $L'$ are homogeneous of degrees $m$, $m'$ then at the
point $P_-$ their expansions (similar to \refE{expan}) begin with
$i=-m-g$, $j=-m'-g$, respectively. The expansion for $\L$ begins
with some integer $m_-$. Hence the condition at $P_-$ looks as
follows:
\[  -m-m'-2g\le \max\,\{0, -1+m_-\}.
\]
Finally, we obtain
\[ \min\,\{0, 1-m_-\}-2g  \le m+m'\le \max\,\{0, -1-m_+\}.
\]
Since $m_\pm$ are fixed ($\L$ is fixed), the latter exactly means
that the cocycle is local.
\end{proof}

\subsection{Central extensions for the Lax operator algebras
     over $\mathfrak{so}(n)$}\label{SS:ext_so2n}

We keep here the same line of arguing as in the previous section.

\begin{lemma}\label{L:eig_so} At every weak singularity, the one-form $\tr\, L\,dL'$ has at most simple pole and
\begin{equation}\label{E:eig}\res\,\tr\, L\,dL'=2\k([L,L']).
\end{equation}
\end{lemma}
\begin{proof}
1) Using \refE{i0'} and the relation
\[ dL'=-L'_{-1}z^{-2}+L'_1+\ldots\ ,
\]
where $L'_{-1}$ is given by \refE{i'}, we obtain
\begin{equation}\label{E:e_so2n}  L\,dL'=-\frac{L_{-1}L'_{-1}}{z^3}- \frac{L_0L'_{-1}}{z^2} - \frac{L_1L'_{-1}-L_{-1}L_1'}{z} + \ldots\ .
\end{equation}
For the first term we have
\[ L_{-1}L_{-1}'=(\a\b^t-\b\a^t)(\a\b'^t-\b'\a^t)
=\a(\b^t\a)\b'^t-\b(\a^t\a)\b'^t-\a\b^t\b'\a^t-\b\a^t\b'\a^t.
\]
The first two summands vanish due to \refE{j'}. For the remainder,
we have
\[ \tr\,(L_{-1}L_{-1}')=\tr\,(-\a\b^t\b'\a^t-\b\a^t\b'\a^t)=
   -(\a^t\a)\b^t\b'-(\a^t\b)\a^t\b'
\]
which again vanishes due to \refE{j'}.

Again (as in \refSS{ext_gln}) the term containing $z^{-2}$ vanishes
under the symbol of trace. By definition,
\[ L_0L'_{-1}=L_0(\a\b'^t-\b'\a^t)=\k\a\b'^t-L_0\b'\a^t.
\]
Now observe that $\tr\,(L_0\b'\a^t)=\tr\,(\a^tL_0\b')$ and
$\a^tL_0=-\k\a^t$. Hence $\tr\,(L_0L'_{-1})=2\k\a^t\b'$ which
vanishes due to \refE{j'}.

The third term in \refE{e_so2n} gives us the required residue. We
have
\[ \res\, L\,dL'= (L_{-1}L_1'-L_1L_{-1}').
\]
Substituting $L_{-1}$, $L_{-1}'$ from \refE{i'}, we have
\[ \res\, L\,dL'= \a\b^tL_1'-\b\a^tL_1'-L_1\a\b'^t+L_1\b'\a^t
\]
hence
\[ \tr\,\res\, L\,dL'= \b^tL_1'\a-\a^tL_1'\b-\b'^tL_1\a+\a^tL_1\b' .
\]
Due to the skew-symmetry of matrices $L_1$, $L_1'$ the first two
summands of the last relation are equal, and the same about the
last two summands. Hence
\[ \tr\,\res\, L\,dL'= 2(\b^tL_1'\a-\b'^tL_1\a) .
\]
From \refE{eigen_so} obviously follows that
$[L,L']_0\a=\b^tL_1'\a-\b'^tL_1\a$ which proves the lemma.
\end{proof}

\begin{lemma}\label{L:iso'}
 Let $\mathcal L$ be a skew-symmetric
 matrix-valued 1-form such that locally, at a weak singularity,
\[
 \L=\L_{-1}\frac{dz}{z}+
    \L_0 dz    +\ldots
\]
where  $\L_{-1}=\a{\tilde\b}^t-{\tilde\b}\a^t$, $\tilde{\b}^t\a=1$
and $\L_0\a={\tilde\k}\a$. Then the $1$-form $\tr\, L\L$ has at
most simple pole at $z=0$ and
\[ \res\,\tr\, L\L =2\k(L).
\]
\end{lemma}
\begin{proof}
The coefficient of the degree $-2$ term in the product $L\L$ is
equal, by definition, to
\[ (\a\b^t-\b\a^t)(\a{\tilde\b}^t-{\tilde\b}\a^t)=
\a(\b^t\a){\tilde\b}^t-\a\b^t{\tilde\b}\a^t-\b(\a^t\a){\tilde\b}^t
 +\b\a^t{\tilde\b}\a^t.
\]
The first and the third terms vanish due to \refE{j'}. For the
trace of the remainder we have
\[ \tr\,(-\a\b^t{\tilde\b}\a^t+\b\a^t{\tilde\b}\a^t)=
  -(\a^t\a)\b^t{\tilde\b}+(\a^t\b)\a^t{\tilde\b}
\]
which vanishes for the same reason.

Multiplying expansions for $L$ and $\L$ we find
\[ \tr\,\res\,(L\L)=\tr\,(L_{-1}\L_0+L_0\L_{-1})= \tr\,(\a\b^t-\b\a^t)\L_0+\tr\, L_0(\a{\tilde\b}^t-{\tilde\b}\a^t).
\]
For the first summand we have
\[ \tr\,(\a\b^t-\b\a^t)\L_0= \tr\,(\a\b^t\L_0-\b\a^t\L_0)=
  \b^t\L_0\a-\a^t\L_0\b.
\]
Observe that, by skew-symmetry, $\a^t\L_0\b=-\b^t\L_0\a$, hence
\[ \tr\,(\a\b^t-\b\a^t)\L_0=\b^t\L_0\a-\a^t\L_0\b= 2\b^t\L_0\a=
 2{\tilde\k}\b^t\a
\]
which vanishes due to \refE{j'}.

For the second summand we have
\[ \tr\, L_0(\a{\tilde\b}^t-{\tilde\b}\a^t)=
  \tr\,( L_0\a{\tilde\b}^t-L_0{\tilde\b}\a^t)=
   {\tilde\b}^tL_0\a-\a^tL_0{\tilde\b}.
\]
Since $L_0\a=\k\a$, $\a^tL_0=-\k\a^t$ and
${\tilde\b}^t\a=\a^t{\tilde\b}=1$ we obtain
\[\tr\, L_0(\a{\tilde\b}^t-{\tilde\b}\a^t)=2\k .
\]
\end{proof}
\begin{theorem}\label{T:loc'}
For every $\L$ satisfying the conditions of the \refL{iso'} the
1-form $\tr\, (L\,dL'-[L,L']\L)$  is regular except at the points
$P_\pm$, and the expression
\[  \ga(L,L')= \res_{P_+}\tr\, (L\,dL'-[L,L']\L)
\]
gives a local cocycle on the Lax operator algebra.
\end{theorem}
The proof is quite similar to the proof of the \refT{loc}. It
relies only on the absence of residues of the $1$-form defining
the cocycle at the weak singularities.

There is a certain ambiguity in the definition of $\L$ in
\refL{iso'}. For example, we could require that
$\L_{-1}=\a{\tilde\b}^t$ and take $L\,dL'-2[L,L']\L$ in the
\refT{loc'}.

\subsection{Central extensions for the Lax operator algebras
     over $\mathfrak{sp}(2n)$}\label{SS:ext_sp2n}

We keep the same line of arguing as in the previous sections.
First, prove the following analogue of the \refL{eig_so}.
\begin{lemma}\label{L:eig_so''} At every weak singularity
the $1$-form $\tr\, L\,dL'$ has at most simple pole, and
\begin{equation}\label{E:eig''}\res\,\tr\, L\,dL'=2\k([L,L']).
\end{equation}
\end{lemma}
\begin{proof} A straightforward calculation of $LdL'$ based on the
expansion \refE{i0''} shows that the terms with $z^{-5}$ and
$z^{-4}$ of that matrix-valued 1-form are equal to $0$ by the
relations \refE{j''}.

For the term with $z^{-3}$ we have
\[ (LdL')_{-3}=-(\b^t\s\b'-2\nu')\a\a^t\s .
\]
This expression vanishes under the symbol of trace:
\[ \tr\,(LdL')_{-3}=-(\b^t\s\b'-2\nu')(\a^t\s\a)=0 .
\]
Similarly, using \refE{j''}-\refE{l''} we have
\[
  \tr\,(LdL')_{-2}=\nu(\a^t\s L_1'\a)-2\k(\b'^t\s\a)-\nu'(\a^t\s L_1\a)=0 .
\]
Thus, the 1-form $\tr\, LdL'$ indeed has at most simple pole at
the point in question. A straightforward calculation results in
\[
   \tr\,(LdL')_{-1}=2(\nu\cdot\a^t\s L_2'\a-\nu'\cdot\a^t\s L_2\a +
   \b^t\s L_1'\a-\b'^t\s L_1\a)
\]
which exactly coincides with the double expression \refE{kappa}
for $\k([L,L'])$.
\end{proof}

\begin{lemma}\label{L:iso''}
 Let $\mathcal L$ be a $\g$-valued 1-form such that locally,
 in the neighborhood of any weak singularity,
\[
 \L=\L_{-1}\frac{dz}{z}+
    \L_0 dz    +\ldots\
\]
where $\L_{-1}=(\a{\tilde\b}^t+{\tilde\b}\a^t)\s$,
$\tilde{\b}^t\s\a=1$, $\L_0\a={\tilde\k}\a$ è $\a^t\s\L_1\a=0$.
Then the $1$-form $\tr\, L\L$ has at most simple pole at $z=0$ è
\[ \res\,\tr\, L\L =2\k(L).
\]
\end{lemma}
\begin{proof}
The expansion for $L\L$ begins with $z^{-3}$. We have
\begin{align*}
  \tr\,(L\L)_{-3} & =\nu({\tilde\b}^t\s\a)(\a^t\s\a)=0,\\
  \tr\,(L\L)_{-2} & =\nu(\a^t\s{\tilde\b})(\a^t\s\a)+
                (\b^t\s{\tilde\b})(\a^t\s\a)+
                (\a^t\s{\tilde\b})(\a^t\s\b)=0.
\end{align*}
Thus, the $1$-form $\tr\, L\L$ indeed has at most simple pole at
the point in question. The calculation of the residue results in
\[
  \tr\,
  (L\L)_{-1}=\nu\cdot\a^t\s\L_1\a+\b^t\s(\L_0\a)+(\a^t\s\L_0)\b+
  {\tilde\b}^t\s(L_0\a)+(\a^t\s L_0){\tilde\b}.
\]
By assumption of the lemma, $\a^t\s\L_1\a=0$,
$\L_0\a={\tilde\k}\a$. The last implies also that
$\a^t\s\L_0=-{\tilde\k}\a^t\s$. For this reason
\[
  \tr\,
  (L\L)_{-1}={\tilde\b}^t\s(L_0\a)+(\a^t\s L_0){\tilde\b}.
\]
By the relations $L_0\a=\k\a$, $\a^t\s L_0=-\k\a^t\s$ we have
\[
  \tr\,(L\L)_{-1}=2\k({\tilde\b}^t\s\a)=2\k.
\]
\end{proof}

As above, the last two lemmas imply
\begin{theorem}\label{T:loc"}
For every $\L$ satisfying the conditions of the \refL{iso''} the
1-form $\tr\, (L\,dL'-[L,L']\L)$  is regular except at the points
$P_\pm$, and the expression
\[  \ga(L,L')= \res_{P_+}\tr\, (L\,dL'-[L,L']\L)
\]
gives a local cocycle on the Lax operator algebra.
\end{theorem}







\end{document}